\documentclass[12pt]{article}
\usepackage{fullpage,graphicx,graphics,amsmath,amsfonts,verbatim,siunitx}

\newcommand{\BEAS}{\begin{eqnarray}}
\newcommand{\EEAS}{\end{eqnarray}}
\newcommand{\BEA}{\begin{eqnarray}}
\newcommand{\EEA}{\end{eqnarray}}
\newcommand{\BEQ}{\begin{equation}}
\newcommand{\EEQ}{\end{equation}}
\newcommand{\BIT}{\begin{itemize}}
\newcommand{\EIT}{\end{itemize}}
\newcommand{\BNUM}{\begin{enumerate}}
\newcommand{\ENUM}{\end{enumerate}}

\newcommand{\BA}{\begin{array}}
\newcommand{\EA}{\end{array}}

\newcommand{\eg}{{\it e.g.}}
\newcommand{\ie}{{\it i.e.}}

\newcommand{\ones}{\mathbf 1}

\newcommand{\reals}{{\mbox{\bf R}}}

\newcommand{\argmin}{\mathop{\rm argmin}}

\newcommand{\intr}{\mathop{\bf int}}

\newcommand{\K}{\mathcal{K}}

\newcounter{algorithmctr}[section]
\renewcommand{\thealgorithmctr}{\thesection.\arabic{algorithmctr}}
\newenvironment{algdesc}%
   {\refstepcounter{algorithmctr}\begin{list}{}{%
       \setlength{\rightmargin}{0\linewidth}%
       \setlength{\leftmargin}{0\linewidth}}%
       \rmfamily\small
       \item[]{\setlength{\parskip}{0ex}\hrulefill\par%
        \nopagebreak{\bfseries\textsf{Algorithm \thealgorithmctr~}}}}%
   {{\setlength{\parskip}{-1ex}\nopagebreak\par\hrulefill} \end{list}}

\newcommand{\eps}{\varepsilon}

\newcommand{\downto}{\downarrow}

\title{Automatic Repair of Convex Optimization Problems}
\author{Shane Barratt \and Guillermo Angeris  \and Stephen Boyd}

\begin{document}
\maketitle

\begin{abstract}
Given an infeasible, unbounded, or pathological convex optimization problem,
a natural question to ask is: what is the smallest change we can make to the problem's parameters such that the problem becomes solvable?
In this paper, we address this question by posing it as an optimization problem
involving the minimization of a convex regularization function of the
parameters, subject to the constraint that the parameters result in a solvable problem.
We propose a heuristic for approximately solving this problem
that is based on the penalty method and leverages recently developed methods
that can efficiently evaluate the derivative of the
solution of a convex cone program with respect to its parameters.
We illustrate our method by applying it to examples in optimal control and economics.
\end{abstract}

\section{Introduction}

\paragraph{Parametrized convex optimization.}
We consider parametrized convex optimization problems, which have the form
\BEQ
\begin{array}{ll}
    \mbox{minimize} & f_0(x;\theta)\\
    \mbox{subject to} & f_i(x;\theta) \leq 0, \quad i=1,\ldots,m, \\
    & g_i(x;\theta)=0,\quad i=1,\ldots,p,
\end{array}
\label{eq:cvxopt}
\EEQ
where $x\in\reals^n$ is the optimization variable, $\theta\in\reals^k$ is the parameter,
the objective function
$f_0:\reals^n\times\reals^k\to\reals$ is convex in $x$,
the inequality constraints functions
$f_i:\reals^n\times\reals^k\to\reals$, $i=1,\ldots,m,$ are
convex in $x$,
and the equality constraint functions $g_i:\reals^n\times\reals^k\to\reals$, $i=1,\ldots,p$,
are affine in $x$.

\paragraph{Solvable problems.}
A point $x$ is said to be feasible if
$f_i(x; \theta) \leq 0$, $i=1,\ldots,m$,
and $g_i(x;\theta)=0$, $i=1,\ldots,k$.
The optimal value $p^\star$ of the problem \eqref{eq:cvxopt} is defined as
\[
p^\star=\inf\{f_0(x;\theta) \mid f_i(x) \leq 0, \; i=1,\ldots,m, \; h_i(x)=0, \; i=1,\ldots,p\}.
\]
We allow $p^\star$ to take on the extended values $\pm \infty$.
Roughly speaking, we say that \eqref{eq:cvxopt} is
\emph{solvable} if $p^\star$ is finite and attainable.
(We will define solvable formally below, when we canonicalize
\eqref{eq:cvxopt} into a cone program.)
When the problem is unsolvable, it falls into one of three cases: it is \emph{infeasible} if $p^\star=+\infty$,
\emph{unbounded below} if $p^\star=-\infty$,
and \emph{pathological} if $p^\star$ is finite,
but not attainable by any $x$, or strong duality does not hold for~\eqref{eq:cvxopt}.
Unsolvable problems
are often undesirable since,
in many cases, there does not exist a solution.

\paragraph{Performance metric.}
The goal in this paper is to repair an unsolvable problem
by adjusting the parameter $\theta$ so that it becomes solvable.
We will judge the desirability of a new parameter $\theta$
by a (convex) performance metric function $r:\reals^k \to \reals \cup \{+\infty\}$,
which we would like to be small.
(Infinite values of $r$ denote constraints on the parameter.)
A simple example of $r$ is the Euclidean distance to an initial parameter vector
$\theta_0$, or $r(\theta)=\|\theta-\theta_0\|_2$.

\paragraph{Repairing a convex optimization problem.}
In this paper, we consider the problem of repairing a convex optimization problem, as measured by the performance metric,
by solving the problem
\BEQ
\begin{array}{ll}
    \mbox{minimize}  & r(\theta)\\
    \mbox{subject to} & \text{problem} \; \eqref{eq:cvxopt} \; \text{is solvable},
\end{array}
\label{eq:goal}
\EEQ
with variable $\theta$.

\paragraph{Pathologies.}
There are various pathologies that can occur in this formulation.
For example, the set of $\theta$ that lead to solvable problems could be open, meaning there
might not exist a solution to \eqref{eq:goal},
or the complement could have (Lebesgue) measure zero,
meaning that the problem can be made solvable
by essentially any perturbation.
Both of these cases can be demonstrated with the following problem:
\BEQ
\begin{array}{ll}
    \mbox{minimize}  & 0\\
    \mbox{subject to} & \theta x = 1,
\end{array}
\label{eq:example}
\EEQ
and regularization function $r(\theta)=\theta^2$.
The set of solvable $\theta$ is $\{\theta \mid \theta \neq 0\}$,
which is both open and has complement with measure zero.
The optimal value of problem \eqref{eq:goal} is 0,
but is not attainable by any solvable $\theta$.
The best we can hope to do in these situations is to produce
a minimizing sequence.

\paragraph{NP-hardness.} Repairing a convex optimization problem
is NP-hard.
To show this, we reduce the 0-1 integer programming problem
\BEQ
\begin{array}{ll}
    \mbox{minimize} & 0\\
    \mbox{subject to} & Ax=b,\\
    & x \in \{0,1\}^n,
\end{array}
\label{eq:zero_one}
\EEQ
with variable $x\in\reals^n$ and data $A\in\reals^{m\times n}$ and $b\in\reals^m$
to an instance of problem~\eqref{eq:goal}.

Let $r(\theta)=0$. The convex optimization problem that we would like to be solvable be
\BEQ
\begin{array}{ll}
    \mbox{minimize} & 0\\
    \mbox{subject to} & x = \theta, \\
    & \theta x = x, \\
    & Ax = b,
\end{array}
\label{eq:np-hard}
\EEQ
with variable $x$.
Problem \eqref{eq:goal} has the same constraints
as \eqref{eq:zero_one}, since \eqref{eq:np-hard} is feasible
if and only if $\theta_i(\theta_i-1)=0$, $i=1,\ldots,n$, and $A\theta = b$.
Therefore, the problem of finding any feasible parameters for~\eqref{eq:cvxopt} (\ie, with $r(\theta) = 0$),
is at least as hard as the 0-1 integer programming problem, which is known to be NP-hard~\cite{karp1972reducibility}.

\section{Cone program formulation}
\label{sec:cone-program}
In practice, most convex optimization problems are solved by reformulating them as equivalent conic programs and passing the numerical data in the reformulated problem to general conic solvers such as SCS~\cite{o2016conic}, Mosek~\cite{mosek}, or Gurobi~\cite{gurobi}. This process of canonicalization is often done automatically by packages like CVXPY~\cite{cvxpy}, which generate a conic program from a high-level description of the problem.

\paragraph{Canonicalization.} For the remainder of the paper, we will consider the canonicalized form of problem~\eqref{eq:cvxopt}.
The primal (P) and dual (D) form of the canonicalized convex cone
program is (see, \eg, \cite{ben2001lectures,boyd2004convex})
\BEQ
\begin{array}{cc}
	\begin{array}{lll}
		\text{(P)} &\mbox{minimize}  & c(\theta)^T x\\
		&\mbox{subject to} &  A(\theta)x+s=b(\theta),\\
		&&s\in  \mathcal{K},
	\end{array}	
	\hspace{1em}&
	\begin{array}{lll}
		\text{(D)}&\mbox{minimize}& -b(\theta)^T y\\
		&\mbox{subject to}& A(\theta)^T y+c(\theta)=0,\\
		&&y\in  \mathcal{K}^*.
	\end{array}
\end{array}
\label{eq:cone_program}
\EEQ
Here $x\in\reals^n$ is the primal variable,
$y\in\reals^m$ is the dual variable,
and $s\in\reals^m$ is the slack variable.
The functions $A:\reals^k \to \reals^{m \times n}$,
$b:\reals^k\to\reals^m$, and $c:\reals^k\to\reals^n$
map the parameter vector $\theta$ in problem \eqref{eq:cvxopt}
to the cone program problem data $A$, $b$, and $c$. 
The set $\mathcal K\subseteq\reals^m$ is a
closed convex cone
with associated dual cone $\mathcal K^*=\{y\mid y^Tx \geq 0 \; \text{for all} \; x \in \mathcal K\}$.

\paragraph{Solution.}
The vector $(x^\star, y^\star, s^\star)$ is a solution
to problem \eqref{eq:cone_program} if
\BEQ
A(\theta)x^\star + s^\star = b(\theta), \quad A(\theta)^Ty^\star + c(\theta) = 0, \quad (s^\star,y^\star) \in \mathcal K\times \mathcal K^*, \quad c(\theta)^Tx^\star + b(\theta)^Ty^\star = 0.
\label{eq:sol}
\EEQ
These conditions merely state that $(x^\star, s^\star)$ is primal feasible,
$y^\star$ is dual feasible, and that there is zero duality gap, which implies that $(x^\star, y^\star, s^\star)$ is optimal by weak duality~\cite[\S5.2.2]{boyd2004convex}.
Problems \eqref{eq:cvxopt} and \eqref{eq:cone_program}
are solvable if and only if there exists a point
that satisfies \eqref{eq:sol}.

\paragraph{Primal-dual embedding.}
The primal-dual embedding of problem \eqref{eq:cone_program}
is the cone program
\BEQ
\begin{array}{ll}
	\mbox{minimize}  & t\\
	\mbox{subject to} & \left\|\begin{bmatrix}A(\theta)x + s - b(\theta) \\ A(\theta)^Ty + c(\theta) \\ c(\theta)^Tx + b(\theta)^Ty\end{bmatrix}\right\|_2 \leq t \\
	& s \in \mathcal K, \quad y \in \mathcal K^*,
\end{array}
\label{eq:repair}
\EEQ
with variables $t$, $x$, $y$, and $s$.
This problem is guaranteed to be feasible
since, for any $\theta \in \reals^k$, setting $x=0$, $y=0$, $s=0$, and $t=\|(b(\theta),c(\theta))\|_2$ yields a feasible point. The problem is also guaranteed to be bounded from below, since the objective is nonnegative.
Taken together, this implies that problem~\eqref{eq:repair} always has a solution, assuming it is not pathological.

\paragraph{Optimal value of \eqref{eq:repair}.}
Let $t^\star:\reals^k \to \reals$
denote the optimal value of problem \eqref{eq:repair}
as a function of $\theta$.
Notably, we have that
$t^\star(\theta) = 0$ if and only if problem \eqref{eq:cvxopt} is solvable,
since if $t^\star(\theta)=0$, the solution
to problem \eqref{eq:repair} satisfies \eqref{eq:sol}
and therefore problem \eqref{eq:cvxopt} is solvable.
On the other hand, if problem \eqref{eq:cvxopt} is solvable,
then there exists a point that
satisfies \eqref{eq:sol} and is feasible for \eqref{eq:repair},
so $t^\star(\theta)=0$.

\paragraph{Differentiability of $t^\star$.}
In practice, $t^\star$ is often a differentiable function of $\theta$.
This is the case when $A$, $b$, and $c$ are differentiable,
which we will assume,
and the optimal value of problem \eqref{eq:repair} is differentiable
in $A$, $b$, and $c$.
Under some technical conditions that are often satisfied in practice,
the optimal value of a cone program
is a differentiable function of its problem data \cite{diffcp2019}.
We will assume that $t^\star$ is differentiable,
and that we can efficiently compute its gradient $\nabla t^\star(\theta)$
using the methods described in \cite{diffcp2019}
and the chain rule.

\paragraph{Reformulation.}
In light of these observations, we can reformulate problem \eqref{eq:goal} as
\BEQ
\begin{array}{ll}
	\mbox{minimize}  & r(\theta)\\
	\mbox{subject to} & t^\star(\theta) = 0,
\end{array}
\label{eq:goal2}
\EEQ
with variable $\theta$.
Here we have replaced the intractable constraint in problem \eqref{eq:goal}
with an equivalent smooth equality constraint.
Since this problem is NP-hard,
we must resort to heuristics to find approximate solutions;
we give one in \S\ref{sec:solution_method}.

\section{Heuristic solution method}
\label{sec:solution_method}

\paragraph{Penalty method.}
One simple heuristic is to use the penalty method to (approximately) solve \eqref{eq:goal2}.
Starting from $\theta^0 \in \reals^k$ and $\lambda^0 > 0$, at iteration $\ell$, the penalty method performs the update 
\BEQ
\theta^{\ell+1} = \underset{\theta}{\mbox{argmin}} \; \lambda^\ell r(\theta) + t^\star(\theta),
\label{eq:update}
\EEQ
and then decreases $\lambda^\ell$, \eg, $\lambda^{\ell+1} = (1/2) \lambda^\ell$, until $t^\star (\theta^\ell) \le \epsilon_\mathrm{out}$ for some given tolerance $\epsilon_\mathrm{out} > 0$.

To perform the update \eqref{eq:update},
we must solve the unconstrained optimization problem
\BEQ
\begin{array}{ll}
	\mbox{minimize} & L(\theta,\lambda^\ell) = \lambda^\ell r(\theta) + t^\star(\theta),
\end{array}
\label{eq:step}
\EEQ
with variable $\theta$.
The objective is the sum
of a differentiable function and a potentially nonsmooth convex function,
for which there exist many efficient methods.
The simplest (and often most effective) of these methods is
the proximal gradient method
(which stems from the proximal
point method \cite{martinet1970breve};
for a modern reference see \cite{nesterov2013gradient}). The proximal gradient method consists of the iterations
\[
\theta^{\ell+1} = \mathbf{prox}_{\alpha^\ell\lambda^\ell r}\left(\theta^\ell - \alpha^\ell
\nabla t^\star(\theta^\ell)\right),
\]
where the proximal operator is defined as
\[
\mathbf{prox}_{\alpha\lambda r}(\tilde \theta) = \argmin_{\theta} \; \alpha\lambda r(\theta) + \frac{1}{2} \|\theta - \tilde \theta\|_2^2.
\]
Since $r$ is convex, evaluating the proximal operator of $\alpha\lambda r$ requires
solving a convex optimization problem.
Indeed, for many practical choices of $r$, its proximal operator
has a closed-form expression \cite[\S6]{parikh2014proximal}.

We run the proximal gradient method until
the stopping criterion
\[
\|(\theta^l - \theta^{l+1})/\alpha^l + (g^{l+1} - g^l)\|_2 \leq \epsilon_\mathrm{in},
\]
is reached,
where $g^l = \nabla t^\star(\theta^l)$,
for some given tolerance $\epsilon_\mathrm{in} > 0$ \cite{barratt2019least}.
We employ the adaptive step size scheme described in \cite{barratt2019least}.
The full procedure is described in algorithm \ref{alg} below.

\begin{algdesc}
\label{alg}
\emph{Finding the closest solvable convex optimization problem.}
\begin{tabbing}
    {\bf given} regularization function $r$, initial parameter $\theta^0$,
    	penalty $\lambda^0$, step size $\alpha^0$, iterations $n_\mathrm{iter}$,\\
    	\qquad outer tolerance $\epsilon_\mathrm{out}$, inner tolerance $\epsilon_\mathrm{in}$.\\
    {\bf for} $l=1,\ldots,n_\mathrm{iter}$ \\
    	\qquad \=\ 1.\ \emph{Compute $t^\star$.} Compute $t^\star(\theta^\ell)$.\\
    	\qquad \=\ 2.\ \emph{Compute gradient of $t^\star$.} Let $g^\ell=\nabla t^\star(\theta^\ell)$.\\
    	\qquad \=\ 3.\ \emph{Compute the gradient step.} Let $\theta^{\ell+1/2}=\theta^\ell - \alpha^\ell g^\ell$. \\
    	\qquad \=\ 4.\ \emph{Compute the proximal operator.} Let $\theta^\mathrm{tent} =\mathbf{prox}_{\alpha^\ell\lambda^\ell r} (\theta^{\ell+1/2})$. \\
    	\qquad \=\ 5.\ {\bf if} $L(\theta^\mathrm{tent}, \lambda^\ell) < L(\theta^\ell,\lambda^\ell)$, \\
    	\qquad \qquad \=\ \emph{Increase step size and accept update}. $\alpha^{\ell+1} = (1.2) \alpha^\ell, \quad \theta^{\ell+1} = \theta^\mathrm{tent}$.\\
    	\qquad \qquad \=\ \emph{Inner stopping criterion}. if $\|(\theta^\ell - \theta^{\ell+1}) / \alpha^\ell + (g^{\ell+1}-g^\ell)\|_2 \leq \epsilon_\mathrm{in}$, then set $\lambda^{\ell+1} = (1/2) \lambda^\ell$. \\
    	\qquad  \=\ 6.\ {\bf else} \emph{Decrease step size and reject update}. $\alpha^{\ell+1} = (1/2) \alpha^\ell, \quad \theta^{\ell+1} = \theta^\ell$. \\
    	\qquad  \=\ 7.\ \emph{Outer stopping criterion}. if $t^\star(\theta^{l+1}) \leq \eps_\mathrm{out}$, quit. \\
    {\bf end for}
\end{tabbing}
\end{algdesc}

\paragraph{Implementation.}
We have implemented algorithm \ref{alg} in Python,
which is available online at
\[
\texttt{https://github.com/cvxgrp/cvxpyrepair}
\]
The interface is the \texttt{repair} method, which,
given a parametrized CVXPY problem \cite{cvxpy}
and a convex regularization function,
uses algorithm \ref{alg} to find the parameters that
approximately minimize that regularization function
and result in a solvable CVXPY problem.
We use SCS~\cite{o2016conic} to solve cone programs
and diffcp~\cite{diffcp2019} to compute the gradient of cone programs.
We require the CVXPY problem to be a disciplined parametrized
program (DPP), so that the mapping from parameters to $(A,b,c)$ is affine,
and hence differentiable
\cite{agrawal2019differentiable}.

Until this point, we have assumed that the optimal value of
problem \eqref{eq:repair} is differentiable in $A$, $b$, and $c$.
However, in our implementation, we do not require this to be the case.
So long as it is differentiable almost everywhere, it is reasonable
to apply the proximal gradient method to \eqref{eq:step}.
At non-differentiable points, we instead compute a heuristic quantity.
For example, a source of non-differentiability is the singularity
of a particular matrix; in this case, diffcp computes
a least-squares approximation of the gradient \cite[\S3]{diffcp2019}.

\section{Examples}

\subsection{Spacecraft landing}
We consider the problem of landing a spacecraft
with a gimbaled thruster.
The dynamics are
\[
	m\ddot x(t) = f(t) - mge_3,
\]
where $m>0$ is the spacecraft mass,
$x(t)\in\reals^3$ is the spacecraft position,
$f(t)\in\reals^3$ is the force applied by the thruster,
$g > 0$ is the gravitational acceleration,
and $e_3=(0,0,1)$.

Our goal, given some initial position $x^\mathrm{init}\in\reals^3$
and velocity $v^\mathrm{init}\in\reals^3$,
is to land the spacecraft at zero position and
velocity at some touchdown time $T > 0$,
\ie, $x(T)=0$ and $\dot x(T)=0$.

We have a total available fuel $M^\mathrm{fuel}$
and a thrust limit $F^\mathrm{max}$.
This results in the constraints
\[
\int_0^T \gamma \|f(t)\|_2\,dt \le M^\mathrm{fuel}, \quad \|f(t)\|_2 \leq F^\mathrm{max}, \quad 0 \le t \le T,
\]
where $\gamma$ is the fuel consumption coefficient.
We also have a gimbal constraint
\[
f_3(t) \geq \alpha \|(f_1(t),f_2(t))\|_2,
\]
where $\alpha$ is equal to the tangent of the maximum
gimbal angle.

We discretize the thrust profile, position, and velocity at intervals of length $h$, or
\[
f_k = f((k-1)h), \quad x_k = x((k-1)h), \quad v_k = \dot x((k-1)h), \quad k=1,\ldots,H,
\]
where $H=T/h + 1$.

To find if there exists a thrust profile to land the spacecraft,
we solve the problem
\BEQ
\begin{array}{ll}
	\mbox{minimize}  & 0\\
	\mbox{subject to} &x_{k+1} = x_k + (h/2)(v_{k+1} + v_{k}), \quad k=1, \dots, H, \\
	& m v_{k+1} = f_k - h m g e_3, \quad k=1, \dots, H, \\
	& \|f_k\|_2 \leq F^\mathrm{max}, \quad k=1,\ldots,H, \\
	& \sum_{k=1}^H h \gamma \|f_k\|_2 \le M^\mathrm{fuel}, \\
	& (f_k)_3 \geq \alpha \|((f_k)_1,(f_k)_2)\|_2, \\
	& x_1 = x^\mathrm{init}, \quad v_1 = v^\mathrm{init},\\
	& x_H = 0, \quad v_H = 0,
\end{array}
\label{eq:spacecraft_problem_discretized}
\EEQ
with variables $x$, $v$, and $f$.
This problem is a parametrized convex optimization problem,
with parameter
\[
\theta = (m, M^\mathrm{fuel}, F^\mathrm{max}, \alpha).
\]

Suppose that we are given a parameter vector $\theta_0$ for which
it is impossible to find a feasible thrust profile, \ie, problem \eqref{eq:spacecraft_problem_discretized}
is infeasible.
Suppose, in addition, that we are allowed to change the spacecraft's parameters
in a limited way.
We seek to find the smallest
changes to the mass and constraints on the fuel and thrust limit
that guarantees the feasibility of problem \eqref{eq:spacecraft_problem_discretized}.
We can (approximately) do this with algorithm \ref{alg}.

\paragraph{Numerical example.} We consider a numerical example with data
\[
T = 10, \quad h = 1, \quad g=9.8, \quad x^\mathrm{init}=(10,10,50), \quad v^\mathrm{init}=(10,-10,-10), \quad \gamma = 1,
\]
and initial parameters
\[
m_0 = 12, \quad M^\mathrm{fuel}_0 = 200, \quad F^\mathrm{max}_0 = 50, \quad \alpha_0=0.5.
\]

The initial parameters are infeasible, \ie, there is no possible thrust profile which allows the spacecraft to land in time,
so we use algorithm~\ref{alg} to modify the design parameters in order to have a feasible landing thrust profile.
We use the performance metric
\[
r(\theta) = \begin{cases}
\frac{|m - m_0|}{m_0} + \frac{|M^\mathrm{fuel} - M^\mathrm{fuel}_0|}{M^\mathrm{fuel}_0} + \frac{|F^\mathrm{max} - F^\mathrm{max}_0|}{F^\mathrm{max}_0} + \frac{|\alpha - \alpha_0|}{\alpha_0} & m \ge 9, \\
+\infty & \text{otherwise},
\end{cases}
\]
which constrains the mass to be greater than or equal to 9,
and penalizes the percentage change in each of the parameters.
The resulting feasible design has the parameters
\[
m=9.03, \quad M^\mathrm{fuel} = 271.35, \quad F^\mathrm{max} = 67.16, \quad \alpha = 0.5,
\]
and $r(\theta) = 0.948$.

\subsection{Arbitrage}
Consider an event (\eg, horse race, sports game, or a financial market over
a short time period)
with $m$ possible outcomes and $n$ possible wagers
on the outcome.
The return matrix is $R\in\reals^{m \times n}$, where
$R_{ij}$ is the return in dollars
for the outcome $i$ and wager $j$ per dollar bet.
A betting strategy is a vector $w\in\reals_+^n$,
where $w_i$ is the amount that we bet on the $i$th wager.
If we use a betting strategy $w$
and outcome $i$ occurs, then the return is $(Rw)_i$ dollars.

We say that there is an \emph{arbitrage opportunity} in this
event if there exists
a betting strategy $w\in\reals_+^n$ that is guaranteed
to have nonnegative return for each outcome, and positive
return in at least one outcome.
We can check whether there exists an arbitrage opportunity
by solving the convex optimization problem
\BEQ
\begin{array}{ll}
	\mbox{maximize}  & \ones^TRw\\
	\mbox{subject to} & Rw \geq 0,\\
    & w \geq 0,
\end{array}
\label{eq:arbitrage}
\EEQ
with variable $w$.
If this problem is unbounded above,
then there is an arbitrage opportunity.

Suppose that we are the event organizer (\eg, sports book director, bookie, or financial exchange)
and we wish to design the return matrix $R$ such that there
is are arbitrage opportunities and that some performance metric
$r$ is small.
We can tackle this problem by finding the
nearest solvable convex optimization problem to problem \eqref{eq:arbitrage}
using algorithm \ref{alg}.

\paragraph{Numerical example.} We consider
a horse race with $n=3$ horses and $m=5$ outcomes.
The initial return matrix is
\[
R_0 = \begin{bmatrix}
  0.05 & 1.74 & -0.88\\
  0.08 & 0.45 & -1.02\\
  0.18 & -0.31 & 1.29\\
  0.9 & -1.17 & 0.27\\
  -0.93 & 0.17 & 2.39\\
\end{bmatrix},
\]
for which there is an arbitrage opportunity in the direction
\[
w = (0.71, 0.62, 0.33).
\]
We consider the regularization function $r(R) = \|(R-R_0)/R_0\|_1$,
where $/$ is meant elementwise.
After running algorithm \ref{alg}, the arbitrage-free return matrix is
\[
R_\mathrm{final} = \begin{bmatrix}
  0.05 & 1.71 & -0.9\\
  0.08 & 0.42 & -1.09\\
  0.18 & -0.31 & 1.27\\
  0.81 & -1.22 & 0.27\\
  -0.97 & 0.17 & 2.37\\
\end{bmatrix},
\quad
R_\mathrm{final} - R_0 = 
\begin{bmatrix}
  0 & -0.03 & -0.02\\
  0 & -0.04 & -0.08\\
  0 & 0 & -0.02\\
  -0.09 & -0.05 & 0\\
  -0.04 & 0 & -0.03\\
\end{bmatrix},
\]
and $r(R_\mathrm{final})=0.142$.

\section{Related work}
Problem~\eqref{eq:goal} is often tractable when, in its conic representation
\eqref{eq:cone_program}, $A$ is constant, and $b$ and $c$ are affine functions
of $\theta$.
(This is not the case in any of our examples.)
In this case, problem~\eqref{eq:goal} can be expressed
as a convex problem.
In the case where the cones are products of the nonnegative reals,
the resulting problem is immediately convex, while the
more general case requires some care (see appendix~\ref{app:conic-solutions}).

This idea is exploited by the Mosek~\cite[\S14.2]{mosek} and
Gurobi~\cite{gurobi} solvers whenever a user would like to
repair an infeasible or unbounded linear program. Automatic
repair for linear programs appears to have been first suggested
in~\cite{roodman1979note}. This was later studied more generally
in the case of linear programs as irreducibly inconsistent systems
(IIS), first defined in~\cite{van1981irreducibly}, with some further
automated repair algorithms in~\cite{chinneck1991locating}, and in
the more general case of linearly-constrained programs in~\cite{leon2001fuzzy}.

The problem that we consider can also be interpreted as automatic
program repair, where the program we are repairing solves
a convex optimization problem \cite{gazzola2017automatic}.
To the best of our knowledge, this paper is the first
to consider automatic convex program repair.

\section*{Acknowledgments}
S. Barratt is supported by the National Science Foundation Graduate Research Fellowship
under Grant No. DGE-1656518.

\bibliography{refs}

\begin{thebibliography}{10}

\bibitem{agrawal2019differentiable}
A.~Agrawal, B.~Amos, S.~Barratt, S.~Boyd, S.~Diamond, and J.~Z. Kolter.
\newblock Differentiable convex optimization layers.
\newblock In {\em Advances in Neural Information Processing Systems}, pages
  9558--9570, 2019.

\bibitem{diffcp2019}
A.~Agrawal, S.~Barratt, S.~Boyd, E.~Busseti, and W.~Moursi.
\newblock Differentiating through a cone program.
\newblock {\em Journal of Applied and Numerical Optimization}, 1(2):107--115,
  2019.

\bibitem{barratt2019least}
S.~Barratt and S.~Boyd.
\newblock Least squares auto-tuning.
\newblock {\em arXiv preprint arXiv:1904.05460}, 2019.

\bibitem{ben2001lectures}
A.~Ben-Tal and A.~Nemirovski.
\newblock {\em Lectures on modern convex optimization: analysis, algorithms,
  and engineering applications}, volume~2.
\newblock SIAM, 2001.

\bibitem{boyd2004convex}
S.~Boyd and L.~Vandenberghe.
\newblock {\em Convex optimization}.
\newblock Cambridge University Press, 2004.

\bibitem{chinneck1991locating}
J.~Chinneck and E.~Dravnieks.
\newblock Locating minimal infeasible constraint sets in linear programs.
\newblock {\em ORSA Journal on Computing}, 3(2):157--168, 1991.

\bibitem{cvxpy}
S.~Diamond and S.~Boyd.
\newblock {CVXPY}: A {P}ython-embedded modeling language for convex
  optimization.
\newblock {\em Journal of Machine Learning Research}, 17(83):1--5, 2016.

\bibitem{gazzola2017automatic}
L.~Gazzola, D.~Micucci, and L.~Mariani.
\newblock Automatic software repair: A survey.
\newblock {\em IEEE Transactions on Software Engineering}, 45(1):34--67, 2017.

\bibitem{gurobi}
{GUROBI Optimization}.
\newblock Gurobi optimizer reference manual.
\newblock 2019.

\bibitem{karp1972reducibility}
R.~Karp.
\newblock Reducibility among combinatorial problems.
\newblock In {\em Complexity of Computer Computations}, pages 85--103.

\bibitem{leon2001fuzzy}
T.~Le{\'o}n and V.~Liern.
\newblock A fuzzy method to repair infeasibility in linearly constrained
  problems.
\newblock {\em Fuzzy Sets and Systems}, 122(2):237--243, 2001.

\bibitem{martinet1970breve}
B.~Martinet.
\newblock Br{\`e}ve communication. r{\'e}gularisation d'in{\'e}quations
  variationnelles par approximations successives.
\newblock {\em ESAIM: Mathematical Modelling and Numerical
  Analysis-Mod{\'e}lisation Math{\'e}matique et Analyse Num{\'e}rique},
  4(R3):154--158, 1970.

\bibitem{mosek}
{MOSEK Aps}.
\newblock {MOSEK} optimizer {API} for {P}ython.
\newblock \texttt{https://docs.mosek.com}, January 2020.

\bibitem{nesterov2013gradient}
Y.~Nesterov.
\newblock Gradient methods for minimizing composite functions.
\newblock {\em Mathematical Programming}, 140(1):125--161, 2013.

\bibitem{o2016conic}
B.~O'Donoghue, E.~Chu, N.~Parikh, and S.~Boyd.
\newblock Conic optimization via operator splitting and homogeneous self-dual
  embedding.
\newblock {\em Journal of Optimization Theory and Applications},
  169(3):1042--1068, 2016.

\bibitem{parikh2014proximal}
N.~Parikh and S.~Boyd.
\newblock Proximal algorithms.
\newblock {\em Foundations and Trends\textregistered{} in Optimization},
  1(3):127--239, 2014.

\bibitem{roodman1979note}
G.~Roodman.
\newblock Note---post-infeasibility analysis in linear programming.
\newblock {\em Management Science}, 25(9):916--922, 1979.

\bibitem{van1981irreducibly}
J.~N.~M. Van~Loon.
\newblock Irreducibly inconsistent systems of linear inequalities.
\newblock {\em European Journal of Operational Research}, 8(3):283--288, 1981.

\end{thebibliography}

\clearpage
\appendix
\section{Convex formulation}
\label{app:conic-solutions}
In the case that $A$ is a constant while $b$ and $c$ are affine functions of $\theta$, we can write~\eqref{eq:goal2} as an equivalent convex optimization problem. In the linear case (\ie, when $\K = \reals_+^n$), we can simply drop the strong duality requirement (which always holds in this case) and express~\eqref{eq:goal2} as
\[
\begin{array}{ll}
	\mbox{minimize}  & r(\theta)\\
	\mbox{subject to} & Ax + s = b(\theta)\\
	& A^Ty + c(\theta) = 0\\
	& s \in \K, \quad y \in \K^*.
\end{array}
\]
For more general cones $\K$ (such as, \eg, the second order cone), a sufficient condition for strong duality is that there exist a feasible point in the interior of the cone. We can write this as, for example, 
\begin{equation}\label{eq:convex-repair}
\begin{array}{ll}
	\mbox{minimize}  & r(\theta)\\
	\mbox{subject to} & Ax + s = b(\theta)\\
	& A^Ty + c(\theta) = 0\\
	& s \in \intr \K, \quad y \in \K^*.
\end{array}
\end{equation}
(We could similarly constrain $y\in\intr\K^*$ and $s\in \K$.)

In general, optimizing over open constraint sets is challenging
and these problems may not even have an optimal point, but,
in practice (and for well-enough behaved $r$, \eg, $r$ continuous)
we can approximate the true optimal value of~\eqref{eq:cone_program}
by approximating the open set $\intr \K$ as a sequence of
closed sets $\K_\eps\subseteq \intr \K$ such
that $\K_\eps \to \intr \K$ as $\eps \downto 0$.

\end{document}